\documentclass[review]{elsarticle}
\usepackage{amsthm}
\usepackage{lineno,hyperref}
\usepackage{amssymb,enumerate}
\usepackage{amsmath}
\usepackage{graphics}
\newtheorem{thm}{Theorem}[section]
 
 \newtheorem{lem}[thm]{Lemma}
 
 \newtheorem{defn}[thm]{Definition}
 \newtheorem{rem}[thm]{Remark}

\journal{Journal of \LaTeX\ Templates}

\begin{document}

\begin{frontmatter}

\title{A note on the blowup of scale invariant damping wave equation with sub-Strauss exponent}

\author{
Ziheng Tu
\footnote{School of Data Science, Zhejiang University of Finance and Economics, 310018, Hangzhou, P.R.China.
 e-mail: tuziheng@zufe.edu.cn.}
\quad
Jiayun Lin
\footnote{
Department of Mathematics and Science, School of Sciences, Zhejiang Sci-Tech University, 310018, Hangzhou, P.R.China.
 email: jylin84@hotmail.com.}
}

\begin{abstract}
We concern the blow up problem to the scale invariant damping wave equations with sub-Strauss exponent. This problem has been studied by Lai, Takamura and Wakasa (\cite{Lai17}) and Ikeda and Sobajima \cite{Ikedapre} recently. In present paper, we extend the blowup exponent from $p_F(n)\leq p<p_S(n+2\mu)$ to $1<p<p_S(n+\mu)$ without small restriction on $\mu$. Moreover, the upper bound of lifespan is derived with uniformly estimate $T(\varepsilon)\leq C\varepsilon^{-2p(p-1)/\gamma(p,n+2\mu)}$. This result extends the blowup result of semilinear wave equation and shows the wave-like behavior of scale invariant damping wave equation's solution even with large $\mu>1$.
\end{abstract}


\end{frontmatter}

\section{Introduction and Main result}
In this paper, we consider the following initial value problem
\begin{equation}\label{main}
\left\{
\begin{array}{ll}
u_{tt}-\Delta u+\frac{\mu}{1+t}u_t=|u|^p\ &(x,t)\ \in\ \mathbb{R}^n\times[0,\infty),\\
u(x,0)=\varepsilon f(x),\ u_t(x,0)=\varepsilon g(x)\ &x\ \in\ \mathbb{R}^n,\\
\end{array}
\right.
\end{equation}
where $\mu>0$, $f,\ g\in C_0^{\infty}(\mathbb{R}^n)$ and $n\in N$. We assume that $\varepsilon>0$ is a "small" parameter.
This type of damping wave equation is called "scale-invariant" due to that the damping term $\frac{\mu}{1+t}u_t$ shares same scaling as $u_{tt}$:
$$\widetilde{u}(t,x)=u(\lambda(1+t)-1,\lambda x).$$
For this typical damping case, the asymptotic behavior of linear equation heavily relies on the size of $\mu$ see \cite{Wirt04}. As far as authors' knowledge, the threshold of $\mu$ according to the asymptotic behavior is still unclear. Meanwhile, the blowup problem or the determination of the critical exponent of the semilinear equation has drawn great of attention. Wakasugi \cite{Wakasugi14} has obtained the blowup result if $1<p\leq p_F(n)$ and $\mu>1$, or $1<p\leq 1+\frac2{n+\mu-1}$ and $0<\mu\leq1$. He has also shown in \cite{WakasugiT} the upper bond of the lifespan:
\begin{equation}\label{Wakasugi}
\left\{
\begin{array}{ll}
C\varepsilon^{-(p-1)/\{2-n(p-1)\}}\ &\mbox{if}\ 1<p<p_F(n)\ \mbox{and}\ \mu>1,\\
C\varepsilon^{-(p-1)/\{2-(n+\mu-1)(p-1)\}}\ &\mbox{if}\ 1<p<1+\frac2{n+\mu-1}\ \mbox{and}\ 0<\mu\leq1,\\
\end{array}
\right.
\end{equation}
where $C$ is a positive constant independent of $\varepsilon$. Here $p_F(n)$ is the Fujita exponent $$p_F(n)=1+\frac2{n}.$$
It is remarkable that, by the so-called Liouville transform:
$$w(x,t):=(1+t)^{\frac{\mu}2}u(x,t),$$
the scale invariant damping wave equation \eqref{main} can be written as Klein Gordon type equation
\begin{equation}\label{Liou}
\left\{
\begin{array}{ll}
w_{tt}-\Delta w+\frac{\mu(2-\mu)}{4(1+t)^2}w=\frac{|w|^p}{(1+t)^{\mu(p-1)/2}} \ &(x,t)\ \in\ \mathbb{R}^n\times[0,\infty),\\
w(x,0)=\varepsilon f(x),\ w_t(x,0)=\varepsilon \{(\mu/2)f(x)+g(x)\}\ &x\ \in\ \mathbb{R}^n.\\
\end{array}
\right.
\end{equation}
Observed that when $\mu=2$, the mass term $\frac{\mu(2-\mu)}{4(1+t)^2}w$ vanishes, so that one can apply some techniques from wave equation. D'Abbicco, Lucente and Reissig \cite{Dabbi15JDE} have obtained following results. Let $\mu=2$, denote the critical exponent
$$p_c(n):=\max\{p_F(n),p_S(n+2)\}$$
where $p_F(n)$ is the Fujita exponent as above and $p_S(n)$ is the Strauss exponent,
$$p_S(n):=\frac{n+1+\sqrt{n^2+10n-7}}{2(n-1)}$$
which is the positive root of the quadratic equation:
$$\gamma(p,n):=2+(n+1)p-(n-1)p^2=0.$$
Then \eqref{main} admits global-in-time solution for sufficiently small $\varepsilon$ if $p>p_c(n)$ in $n=2,\ 3$ though radial symmetry is required in case $n=3$. Hence combing the blowup result from Wakasugi \cite{Wakasugi14}, for this typical case $\mu=2$, dimension $n=2$, the critical exponent is determined. In case of dimension $n=1$ and $\mu=2$, Wakasa \cite{Wakasa16} has verified the critical exponent $p_c(1)=p_F(1)=3$ and showed the optimal of lifespan. Besides, he also showed the critical exponent changes to $p_S(1+2)$ when the nonlinearity is a sign changing type as $|u|^{p-1}u$ and the initial data is of odd function.

Recently Lai, Takamura and Wakasa \cite{Lai17} found such Strauss type exponent exists not only for this specific case $\mu=2$ but also for $\mu$ in range $(0,\frac{n^2+n+2}{2(n+2)})$. In fact, they obtained following result on the blowup exponent and the lifespan:
$$\mbox{for}\ \ \  p_F(n)\leq p< p_S(n+2\mu)\ \ \mbox{and}\ \ 0<\mu<\frac{n^2+n+2}{2(n+2)},$$
$$\mbox{with lifespan}\ \ \ T(\varepsilon)\leq C\varepsilon^{-2p(p-1)/\gamma(p,n+2\mu)}.$$
This exploring shows the wave like behavior appears even for large quantity of $\mu>1$ concerning its blowup phenomena. Very recently, Ikeda and Sobajima \cite{Ikedapre} extended this result to:
\begin{equation*}
T(\varepsilon)\leq\left\{
\begin{array}{ll}
\exp(C\varepsilon^{-p(p-1)})\ \ \ \ \ \ \ \ \ \ &\mbox{for}\ p=p_S(n+\mu)\\
C_{\delta} \varepsilon^{-2p(p-1)/\gamma(p,n+\mu)-\delta}\ \ &\mbox{for}\ p_S(n+2+\mu)\leq p<p_S(n+\mu)\\
C_{\delta}' \varepsilon^{-1-\delta}\ \ \ \ \ \ \ \ \ \ \  \ \ \ \ \ \ &\mbox{for}\ p_F(n)\leq p<p_S(n+2+\mu)\\
\end{array}\right.
\end{equation*}
when $n\geq 2$ and $0\leq \mu<\mu_*:=\frac{n^2+n+2}{n+2}$ and
\begin{equation*}
T(\varepsilon)\leq\left\{
\begin{array}{ll}
\exp(C\varepsilon^{-p(p-1)})\ \ \ \ &\mbox{for}\ p=p_S(1+\mu)\\
C_{\delta} \varepsilon^{-2p(p-1)/\gamma(p,1+\mu)-\delta}\ \ \ \ &\mbox{for}\ \max\{3,\frac2\mu\}\leq p<p_S(1+\mu)\\
C_{\delta}' \varepsilon^{-2(p-1)/\mu-\delta}\ \ \ \ &\mbox{for}\ 0<\mu<\frac23,\ 3\leq p<\frac2\mu\\
\end{array}
\right.
\end{equation*}
when $n=1$ and $0<\mu<\frac{4}{3}$, with arbitrary small $\delta>0$. Their proof relies on the use of hypergeometric function, which is initiated from Zhou-Han \cite{Zhou14}. Their proof deals with critical and sub-critical Strauss exponent cases in a unified way which is quite concise.

In present paper, we consider this blowup problem again. By applying test function method and iteration argument, we improve the above results. Our main novelty is to introduce the modified Bessel function of second kind $K_\nu(z)$. This idea comes from the study of blowup problem of Tricomi equation. He, Witt and Yin \cite{Tri1} used such type special function as test function to derive the blowup exponent of generalized Tricomi equation:
\begin{equation*}
\left\{
\begin{array}{ll}
u_{tt}-t^m\Delta u=|u|^p\ &(x,t)\in\ \mathbb{R}^n\times[0,\infty),\\
u(x,0)=\varepsilon f(x),\ u_t(x,0)=\varepsilon g(x)\ &x\in\ \mathbb{R}^n.\\
\end{array}
\right.
\end{equation*}
Inspiring by this, the function $\lambda(t):=(1+t)^{\frac{\mu+1}2}K_{\frac{\mu-1}2}(1+t)$ is found and the test function is constructed which in turn satisfies the conjugate equation of scale invariant damping wave equation
$$u_{tt}-\Delta u-(\frac{\mu}{1+t}u)_t=0.$$
Consequently, a better lower bound estimate of related functional is obtained. We emphasis this estimation is crucial to extending the blowup exponent range. For the proof of main theorem, we follows the iteration arguments in \cite{Lai-scattering} where Lai and Takamura showed the blowup for the scattering damping wave equation with sub-Strauss exponent.

We now state the definition of energy solution and the main result.
\begin{defn}
We say that $u$ is an energy solution of \eqref{main} on $[0,T)$ if
$$u\in C([0,T),H^1(\mathbb{R}^n))\cap C^1([0,T),L^2(\mathbb{R}^n))\cap L^p_{loc}(\mathbb{R}^n\times[0,T))$$
and satisfies
\begin{eqnarray}
&&\int_{\mathbb{R}^n}u_t(x,t)\phi(x,t)dx-\int_{\mathbb{R}^n}u_t(x,0)\phi(x,0)dx\nonumber\\
&&+\int_0^tds\int_{\mathbb{R}^n}\{-u_t(x,s)\phi_t(x,s)+\nabla u(x,s)\cdot\nabla\phi(x,s)\}dx\label{def}\\
&&+\int_0^t\int_{\mathbb{R}^n}\frac{\mu u_t(x,s)}{1+s}\phi(x,s)dx=\int_0^tds\int_{\mathbb{R}^n}|u(x,s)|^p\phi(x,s)dx\nonumber
\end{eqnarray}
with any $\phi\in C_0^\infty(\mathbb{R}^n\times[0,T))$ and any $t\in [0,T)$.
\end{defn}
By employing integration by parts in \eqref{def} and letting $t\rightarrow T$, we have exactly the definition of a weak solution of \eqref{main}.
Our main result is stated in the following.
\begin{thm}
Let $n\geq2,\ \mu>0$ and $1<p<p_S(n+\mu)$. Assume that both $f\in\ H^1(\mathbb{R}^n)$ and $g\in L^2(\mathbb{R}^n)$ are nonnegative and do not vanish identically. Suppose that an energy solution $u$ of \eqref{main} satisfies
$$\mbox{\rm supp}\ u\ \subset\{(x,t)\in \mathbb{R}^n\times [0,T): |x|\leq t+R \}  $$
with some $R\geq1$. Then there exists a constant $\varepsilon_0=\varepsilon_0(f,g,n,p,\mu,R)>0$ such that $T$ has to satisfy
$$T\leq C\varepsilon ^{-2p(p-1)/\gamma(p,n+\mu)}$$
for $0<\varepsilon\leq\varepsilon_0$, where $C$ is a positive constant independent of $\varepsilon$.
\end{thm}
Our result improved Ikeda and Sobajima's results in several ways. We removed the arbitrary small $\delta$ away and in the range $p_F(n)\leq p<p_S(n+2+\mu)$, the lifespan we provided is better, since
$$\gamma(p,n+2+\mu)>0\Leftrightarrow -1<-\frac{2p(p-1)}{\gamma(p,n+\mu)}.$$
Besides, in our results, there is no small restriction on $\mu$ and the lower range of $p$ can be extended to $1$.
\begin{rem}
Our results cover the super-Fujita range in Wakasugi \cite{WakasugiT}. Moreover, the lifespan estimates are also updated in some sub-Fujita range for $\mu<\mu_*$. Specifically, we have following observation.\\
For $1<\mu<\mu_*$ and $p\in (1,p_F(n))$, assume that
\begin{equation}\frac{2p(p-1)}{\gamma(p,n+\mu)}<\frac{p-1}{2-n(p-1)}\label{ass}\end{equation}
which implies $p>\frac2{n+1-\mu}$. Combining $p<p_F(n)$, it is necessary to require
$$1+\frac2n>\frac2{n+1-\mu}$$
which is automatically satisfied by $\mu<\mu_*$. Hence, the assumption \eqref{ass} always holds in the case $1<\mu<\mu_*$ and $p\in(\max (1,\frac2{n+1-\mu}),p_F(n))$.\\
For $0<\mu\leq1$ and $p\in (1,1+\frac2{n+\mu-1})$, we assume
\begin{equation}\label{ass2}\frac{2p(p-1)}{\gamma(p,n+\mu)}<\frac{p-1}{2-(n+\mu-1)(p-1)}\end{equation}
which implies $p>\frac2{n+\mu-1}$. Hence, for $p\in (\max(1,\frac2{n+\mu-1}), 1+\frac2{n+\mu-1})$, the assumption \eqref{ass2} holds.
\end{rem}
As the blowup result of Strauss critical exponent $p=p_S(n+\mu)$ has been given by Ikeda and Sobajima \cite{Ikedapre}, we note that concerning the determining of critical exponent of scale invariant damping wave equation, the situation of $p>p_s(n+\mu)$ for $0<\mu\leq\mu^*$ needs to be further investigated.

\section{Preliminaries}
Let $u$ be an energy solution of \eqref{main} on $[0,T)$ and define the functional
$$G(t):=\int_{\mathbb{R}^n}u(x,t)dx.$$
Choosing the test function $\phi=\phi(x,s)$ in \eqref{def} to satisfy $\phi\equiv1$ in $\{(x,s)\in \mathbb{R}^n \times [0,t]:|x|\leq s+R\}$, we obtain
$$\int_{\mathbb{R}^n}u_t(x,t)dx-\int_{\mathbb{R}^n}u_t(x,0)dx+\int_0^tds\int_{\mathbb{R}^n}\frac{\mu u_t(x,s)}{1+s}dx=\int_0^tds\int_{\mathbb{R}^n}|u(x,s)|^pdx$$
which means that
$$G'(t)-G'(0)+\int_0^t\frac{\mu G'(s)}{1+s}ds=\int_0^tds\int_{\mathbb{R}^n}|u(x,s)|^pdx.$$
Since all the quantities in this equation except $G'(t)$ is differentiable in $t$, so that so is $G'(t)$. Hence, we have
\begin{equation*}
G''(t)+\frac{\mu}{1+t}G'(t)=\int_{\mathbb{R}^n}|u(x,t)|^pdx.
\end{equation*}
Multiplying $(1+t)^\mu$ and then integrating over $[0,t]$, we arrive at the identity
\begin{equation} \label{key}
(1+t)^\mu G'(t)-G'(0)=\int_0^t(1+s)^\mu ds\int_{\mathbb{R}^n}|u(x,s)|^pdx.
\end{equation}
By the positivity assumption on initial data, further integration on $[0,t]$ gives
\begin{eqnarray}
 G(t)&\geq&\int_0^t(1+\tau)^{-\mu}d\tau\int_0^\tau (1+s)^\mu ds\int_{\mathbb{R}^n}|u(x,s)|^pdx\label{iter1}\\
  &\geq&C_0\int_0^td\tau\int_0^\tau (\frac{1+s}{1+\tau})^{\mu}(1+s)^{n(1-p)}|G(s)|^pds\label{iter2}
\end{eqnarray}
where the H\"{o}lder inequality and compact support of solution is used in second line and
$$C_0:=\{\mbox{vol}(B^n(0,1))\}^{1-p}R^{-n(p-1)}>0.$$
In order to initiate the iteration procedure, we also need to give the low bound of $\int_{\mathbb{R}^n}|u(x,t)|^pdx$ in \eqref{iter1}.
In fact, we have following lemma.
\begin{lem}
Suppose the Cauchy problem \eqref{main} has an energy solution $u$ with the initial data $f$ and $g$ satisfying the assumption of Theorem 1.2,
then there exists large $T_0$ which is independent with $f,\ g$ and $\varepsilon$, for any $t>T_0$ and $p>1$,
\begin{equation}\label{Priori}
\int_{\mathbb{R}^n}|u(x,t)|^pdx\geq C_1\varepsilon^p(1+t)^{n-1-\frac{n+\mu-1}2 p}
\end{equation}
where $C_1:=\frac12C_{f,g}^{p}C_{\varphi,R}^{1-p}e^{p(1-R)}\pi^{-p}$.
\end{lem}
Before give the proof of this lemma, we would first introduce the test function. Let $K_{\nu}(t)$ be the modified Bessel function of second kind
$$K_{\nu}(t)=\int_0^{\infty}\exp(-t\cosh z)\cosh(\nu z)dz,\ \ \nu\in R,$$
which is a solution of the equation
$$\bigg(t^2\frac{d^2}{dt^2}+t\frac{d}{dt}-(t^2+\nu^2)\bigg)K_{\nu}(t)=0, \ \ t>0.$$
From \cite{Erdelyi}, page 24, we have
\begin{equation}K_{\nu}(t)=\sqrt{\frac{\pi}{2t}}e^{-t}(1+O(t^{-1}))\ \ \ \mbox{as}\ t\rightarrow\infty.\label{K_v}\end{equation}
Moreover, its derivative identity holds:
\begin{eqnarray}\label{K3}
\frac{d}{dt}K_{\nu}(t)&=&-K_{\nu+1}(t)+\frac{\nu}{t}K_{\nu}(t)\\
&=&-\frac12(K_{\nu+1}(t)+K_{\nu-1}(t)).\label{K4}
\end{eqnarray}
Now we set
$$\lambda(t)=(1+t)^{\frac{\mu+1}{2}}K_{\frac{\mu-1}2}(1+t).$$
It is clear by direct computation that $\lambda(t)$ satisfies
\begin{equation}
\begin{cases}
\begin{aligned}\label{lamb}
&\bigg((1+t)^2\frac{d^2}{dt^2}-\mu(1+t)\frac{d}{dt}+(\mu-(1+t)^2)\bigg)\lambda(t)=0, \ \ t>0.\\
&\lambda(0)=K_{\frac{\mu-1}2}(1),\ \ \ \lambda(\infty)=0.
\end{aligned}
\end{cases}
\end{equation}
Let
$$\varphi(x)=\int_{\mathbb{S}^{n-1}}e^{x\cdot \omega}d\omega,$$
where $\varphi(x)$ satisfies
\begin{equation}
\varphi(x)\sim C_n|x|^{-\frac{n-1}2}e^{|x|} \ \ \ \mbox{as}\ \ \ |x|\rightarrow\infty.
\end{equation}
Also, it is known
$$\Delta\varphi(x)=\varphi(x).$$
We then define the test function
$$\psi(t,x)=\lambda(t)\varphi(x).$$
Now we can give the proof of Lemma 2.1.
\begin{proof}
Define the functional
$$G_1(t):=\int_{\mathbb{R}^n}u(x,t)\psi(x,t)dx$$
with $\psi(t,x)$ defined above. Then by H\"{o}lder inequality, we have
\begin{equation}\int_{\mathbb{R}^n}|u(x,t)|^pdx\geq\frac{|G_1(t)|^p}{(\int_{|x|\leq t+R}\psi^{p'}(t,x)dx)^{p-1}}.\label{factor}\end{equation}
Following we estimate the lower bound of $|G_1(t)|$ and upper bound of $\int_{|x|\leq t+R}\psi^{p'}(t,x)dx$ respectively.
From the definition of energy solution, we have
\begin{eqnarray*}
&&\int_0^t\int_{\mathbb{R}^n}u_{tt}\psi dxds-\int_{0}^t\int_{\mathbb{R}^n}u\Delta\psi dxds\\
&&+\int_0^t\int_{\mathbb{R}^n}\partial_s(\frac{\mu}{1+s}\psi u)-\partial_s(\frac{\mu}{1+s}\psi)udxds=\int_0^t\int_{\mathbb{R}^n}|u|^p\psi dxds.
\end{eqnarray*}
Applying the integration by parts and $\Delta\varphi(x)=\varphi$, we obtain:
\begin{eqnarray*}&&\int_0^t\int_{\mathbb{R}^n}u_{tt}\psi dxds+\int_{0}^t\int_{\mathbb{R}^n}u\varphi(-\lambda+\frac{\mu}{(1+s)^2}\lambda-\frac{\mu}{1+s}\lambda')dxds\\
&&+\int_{\mathbb{R}^n}\frac{\mu}{1+s}\psi udx\bigg|_0^t =\int_0^t\int_{\mathbb{R}^n}|u|^p\psi dxds.
\end{eqnarray*}
Due to \eqref{lamb}, the above equation simplifies to
$$\int_0^t\int_{\mathbb{R}^n}u_{tt}\psi dxds-\int_{0}^t\int_{\mathbb{R}^n}u\varphi\lambda''dxds+\int_{\mathbb{R}^n}\frac{\mu}{1+s}\psi udx\bigg|_0^t =\int_0^t\int_{\mathbb{R}^n}|u|^p\psi dxds.$$
Thus the integration by parts gives
$$\int_{\mathbb{R}^n}(u_t\psi-u\psi_t+\frac{\mu}{1+s}u\psi) dx\bigg|_0^t=\int_0^t\int_{\mathbb{R}^n}|u|^p\psi dxds.$$
As the righthand side integral is positive, we obtain
$$G_1'(t)+\big(\frac{\mu}{1+t}-2\frac{\lambda'(t)}{\lambda(t)}\big)G_1(t)\geq\varepsilon \int_{\mathbb{R}^n}\bigg(g(x)\lambda(0)+(\mu\lambda(0)-\lambda'(0))f(x)\bigg)\varphi(x)dx.$$
Since
\begin{eqnarray*}
\lambda'(t)&=&\frac{\mu+1}2(1+t)^{\frac{\mu-1}2}K_{\frac{\mu-1}2}(1+t)+(1+t)^{\frac{\mu+1}2}K'_{\frac{\mu-1}2}(1+t)\\
&=&\mu(1+t)^{\frac{\mu-1}2}K_{\frac{\mu-1}2}(1+t)-(1+t)^{\frac{\mu+1}2}K_{\frac{\mu+1}2}(1+t),
\end{eqnarray*}
we obtain
$$\mu\lambda(0)-\lambda'(0)=\mu\lambda(0)-(\mu\lambda(0)-K_{\frac{\mu+1}2}(1))=K_{\frac{\mu+1}{2}}(1)>0.$$
Denote $$C_{f,g}:=\int_{\mathbb{R}^n}\bigg(g(x)\lambda(0)+(\mu\lambda(0)-\lambda'(0))f(x)\bigg)\varphi(x)dx,$$
then by the compact support of $g(x)$ and $f(x)$, $C_{f,g}$ is finite and positive.
We come to the differential inequality of $G_1$
$$G_1'(t)+\big(\frac{\mu}{1+t}-2\frac{\lambda'(t)}{\lambda(t)}\big)G_1(t)\geq\varepsilon C_{f,g}.$$
Multiplying $\frac{(1+t)^\mu}{\lambda^2(t)}$ on two sides and then integrating over $[0,t]$, we derive
$$G_1(t)\geq\varepsilon C_{f,g}\frac{\lambda^2(t)}{(1+t)^\mu}\int_0^t\frac{(1+s)^\mu}{\lambda^2(s)}ds.$$
Inserting $\lambda(t)=(1+t)^{\frac{\mu+1}{2}}K_{\frac{\mu-1}2}(1+t)$, we obtain the lower bound of $G_1$
\begin{equation} G_1(t)\geq\varepsilon C_{f,g}\int_0^t\frac{(1+t)K^2_{\frac{\mu-1}{2}}(1+t)}{(1+s)K^2_{\frac{\mu-1}{2}}(1+s)}ds.\label{G1}\end{equation}
The denominator of \eqref{factor} can be estimated in standard way.
\begin{eqnarray}
&&\int_{|x|\leq t+R}\psi^{p'}(t,x)dx\leq \lambda^\frac{p}{p-1}(t)\int_{|x|\leq t+R}\varphi^{p'}(x)dx\nonumber\\
&&\leq (1+t)^{\frac{\mu+1}{2}\frac{p}{p-1}}K^{\frac{p}{p-1}}_{\frac{\mu-1}2}(1+t)\cdot C_{\varphi}(R+t)^{n-1-\frac{n-1}2\frac{p}{p-1}}e^{\frac{p}{p-1}(t+R)}\nonumber\\
&&\leq C_{\varphi,R}(1+t)^{n-1+\frac{p(\mu+1)-(n-1)p}{2(p-1)}}e^{\frac{p}{p-1}(t+R)}K^{\frac{p}{p-1}}_{\frac{\mu-1}2}(1+t),\label{deno}
\end{eqnarray}
where $C_{\varphi,R}=\max(C_{\varphi} R^{n-1-\frac{(n-1)p}{2 (p-1)}},\ C_{\varphi})$.\\
Combing the estimate \eqref{G1}, \eqref{deno} and \eqref{factor}, we now have
\begin{eqnarray*}&&\int_{\mathbb{R}^n}|u(x,t)|^pdx\\
&&\geq\frac{\varepsilon^p C^p_{f,g}(1+t)^pK^{2p}_{\frac{\mu-1}{2}}(1+t)(\int_0^t\frac{1}{(1+s)K^2_{\frac{\mu-1}{2}}(1+s)}ds)^p}
{C^{p-1}_{\varphi,R}(1+t)^{(n-1)(p-1)+\frac{p(\mu+1)-(n-1)p}{2}}e^{p(t+R)}K^{p}_{\frac{\mu-1}2}(1+t)}\\
&&\geq C_{f,g}^{p}C_{\varphi,R}^{1-p}e^{p(1-R)}\varepsilon^p(1+t)^{\frac p2(2-n-\mu)+(n-1)}\\
&& \cdot e^{-p(t+1)}K^p_{\frac{\mu-1}2}(1+t)(\int_0^t\frac1{(1+s)K^2_{\frac{\mu-1}2}(1+s)}ds)^p.
\end{eqnarray*}
Since \eqref{K_v}, then for sufficient large $T_0(>2)$ (which is independent with $f,\ g,\ \varepsilon$) and $t>T_0$, we have
$$K^p_{\frac{\mu-1}2}(1+t)\sim (\frac{\pi}{2(1+t)})^\frac p2 e^{-p(t+1)}$$
and
\begin{eqnarray*}&&(\int_0^t\frac1{(1+s)K^2_{\frac{\mu-1}2}(1+s)}ds)^p\geq (\int_{\frac t2}^t\frac{2}\pi e^{2(1+s)}ds)^p\\
&&=[\frac1\pi(e^{2(1+t)}-e^{2+t})]^p\geq\frac1{2\pi^p}e^{2p(1+t)}.
\end{eqnarray*}
Consequently,
$$\int_{\mathbb{R}^n}|u(x,t)|^pdx\geq C_{1}\varepsilon^p(1+t)^{\frac p2(1-n-\mu)+(n-1)}\ \ \mbox{for}\ \ t>T_0,$$
where $$C_1:=\frac12C_{f,g}^{p}C_{\varphi,R}^{1-p}e^{p(1-R)}\pi^{-p}.$$
\end{proof}

\section{Proof of Main Theorem}
In this section, we devote to prove Theorem 1.2. The iteration method is applied based on the low bound estimate \eqref{iter1}, \eqref{iter2} and Lemma 2.1.
\begin{proof}
Plugging \eqref{Priori} into \eqref{iter1}, we have for $t>T_0$,
\begin{eqnarray*}
G(t)&\geq&\int_0^t(1+\tau)^{-\mu}d\tau\int_0^\tau(1+s)^\mu C_1\varepsilon^p(1+s)^{n-1-\frac{n+\mu-1}{2}p}ds\\
&\geq&C_1\varepsilon^p\int_{T_0}^t(1+\tau)^{-\mu}d\tau\int_{T_0}^\tau(1+s)^{(n+\mu-1)(1-\frac p2)}ds\\
&\geq&C_1\varepsilon^p\int_{T_0}^t(1+\tau)^{-\mu-(n+\mu-1)\frac p2}d\tau\int_{T_0}^\tau(1+s)^{n+\mu-1}ds\\
&\geq&C_1\varepsilon^p(1+t)^{-\mu-(n+\mu-1)\frac p2}\int_{T_0}^td\tau\int_{T_0}^\tau(s-T_0)^{n+\mu-1}ds.
\end{eqnarray*}
That is
\begin{equation}\label{j=1}
G(t)\geq C_2\varepsilon^p(1+t)^{-\mu-(n+\mu-1)\frac p2}(t-T_0)^{n+\mu+1}\ \ \mbox{for}\ \ t>T_0
\end{equation}
where $C_2=\frac{C_1}{(n+\mu)(n+\mu+1)}$.

Now we begin our iteration argument. Assume that
\begin{equation}\label{iter-assu}
G(t)>D_j(1+t)^{-a_j}(t-T_0)^{b_j}\ \ \ \mbox{for}\ \ t>T_0,\ \ j=1,2,3\ \cdots
\end{equation}
with positive constants $D_j,\ a_j$ and $b_j$ determined later. \eqref{j=1} asserts \eqref{iter-assu} is true for $j=1$ with
\begin{equation}\label{series1}D_1=C_2\varepsilon^p,\ \ a_1=\mu+(n+\mu-1)\frac p2,\ \ \ b_1=n+\mu+1.\end{equation}
Plugging \eqref{iter-assu} into \eqref{iter2}, we have for $t>T_0$
\begin{eqnarray*}
G(t)&>&C_0\int_0^t(1+\tau)^{-\mu}d\tau\int_0^\tau(1+s)^{\mu+n(1-p)}|G(s)|^pds\\
&>&C_0\int_{T_0}^t(1+\tau)^{-\mu}d\tau\int_{T_0}^\tau(1+s)^{\mu+n(1-p)}D_j^p(1+s)^{-pa_j}(s-T_0)^{pb_j}ds\\
&>&C_0D_j^p(1+t)^{-\mu-n(p-1)-pa_j}\int_{T_0}^t(s-T_0)^{\mu+pb_j}dsd\tau\\
&>&\frac{C_0D_j^p}{(\mu+pb_j+1)(\mu+pb_j+2)}(1+t)^{-\mu-n(p-1)-pa_j}(t-T_0)^{\mu+pb_j+2}.
\end{eqnarray*}
So the assumption \eqref{iter-assu} is true if the sequence $\{D_j\}$, $\{a_j\}$, $\{b_j\}$ are define by
\begin{equation}\label{series2}D_{j+1}\geq\frac{C_0}{(\mu+pb_j+2)^2}D_j^p,\ \ a_{j+1}=\mu+n(p-1)+pa_j,\ \ b_{j+1}=\mu+2+pb_j.\end{equation}
It follows from \eqref{series1} and \eqref{series2} that for $j=1,2,3\cdots$
\begin{eqnarray}
a_j&=&[\mu+(n+\mu-1)\frac p2+n+\frac{\mu}{p-1}]p^{j-1}-(n+\frac{\mu}{p-1})\nonumber\\
&=&\alpha p^{j-1}-(n+\frac{\mu}{p-1})\label{a_j}\\
b_j&=&[n+\mu+1+\frac{\mu+2}{p-1}]p^{j-1}-\frac{\mu+2}{p-1}\nonumber\\
&=&\beta p^{j-1}-\frac{\mu+2}{p-1}\label{b_j}
\end{eqnarray}
where we denote the positive constants $$\alpha=\mu+(n+\mu-1)\frac p2+n+\frac{\mu}{p-1},$$ $$\beta=n+\mu+1+\frac{\mu+2}{p-1}.$$
We employ the inequality
$$b_{j+1}=pb_j+\mu+2<p^j[n+\mu+1+\frac{\mu+2}{p-1}],$$
for $D_{j+1}$ to obtain
$$D_{j+1}\geq C_3\frac{D^p_j}{p^{2j}}$$
where $$C_3=\frac{C_0}{(n+\mu+1+\frac{\mu+2}{p-1})^2}.$$
Hence,
\begin{eqnarray*}
\log D_j&\geq& p\log D_{j-1}-2(j-1)\log p+\log C_3\\
&\geq& p^2\log D_{j-2}-2(p(j-2)+(j-1))\log p+(p+1)\log C_3\\
&\geq&\cdots\\
&\geq&p^{j-1}\log D_1-2\log p\sum_{k=1}^{j-1}kp^{j-1-k}+\log C_3\sum_{k=1}^{j-1}p^k.
\end{eqnarray*}
Direct calculation gives
$$\sum_{k=1}^{j-1}kp^{j-1-k}=\frac{1}{p-1}(\frac{p^j-1}{p-1}-j)$$ and
$$\sum_{k=1}^{j-1}p^k=\frac{p-p^j}{1-p},$$
which yields
\begin{eqnarray*}
\log D_j&\geq& p^{j-1}\log D_1-\frac{2\log p}{p-1}(\frac{p^j-1}{p-1}-j)+\log C_3\frac{p-p^j}{1-p}\\
&=&p^{j-1}\bigg(\log D_1-\frac{2p\log p}{(p-1)^2}+\frac{p\log C_3}{p-1}\bigg)+\frac{2\log p}{p-1}j+\frac{2\log p}{(p-1)^2}+\frac{p\log C_3}{1-p}
\end{eqnarray*}
Consequently for $j>\left[\frac{p\log C_3}{2\log p}-\frac{1}{p-1}\right]+1$,
\begin{equation}\label{D_j}D_j\geq\exp\{p^{j-1}(\log D_1-S_p(\infty))\}\end{equation}
with
$$ S_p(\infty):=\frac{2p\log p}{(p-1)^2}-\frac{p\log C_3}{p-1}.$$
Inserting \eqref{a_j}, \eqref{b_j} and \eqref{D_j} into \eqref{iter-assu} gives
\begin{eqnarray}
G(t)&\geq&\exp\big(p^{j-1}(\log D_1-S_p(\infty))\big)(1+t)^{-\alpha p^{j-1}+(n+\frac{\mu}{p-1})}(t-T_0)^{\beta p^{j-1}-\frac{\mu+2}{p-1}}\nonumber\\
&\geq&\exp\big(p^{j-1}J(t)\big)(1+t)^{n+\frac\mu{p-1}}(t-T_0)^{-\frac{\mu+2}{p-1}}\label{contr}
\end{eqnarray}
where
$$J(t):=\log D_1-S_p(\infty)-\alpha\log(1+t)+\beta\log(t-T_0).$$
For $t>2T_0+1$, we have
\begin{eqnarray*}
J(t)&\geq& \log D_1-S_p(\infty)-\alpha\log(2t-2T_0)+\beta\log(t-T_0)\\
&\geq&\log D_1-S_p(\infty)+(\beta-\alpha)\log(t-T_0)-\alpha\log2\\
&=&\log(D_1\cdot(t-T_0)^{\beta-\alpha})-S_p(\infty)-\alpha\log2.
\end{eqnarray*}
Note that
$$\beta-\alpha=\frac{p+1}{p-1}-(n+\mu-1)\frac p2=\frac{\gamma(p,n+\mu)}{2(p-1)}.$$
Thus if $$t>\max\{T_0+(\frac{e^{[S_p(\infty)+\alpha\log2]+1}}{C_2\varepsilon^p})^{2(p-1)/\gamma(p,n+\mu)},2T_0+1\},$$
we then get $J(t)>1$, and this in turn give $G(t)\rightarrow\infty$ by taking $j\rightarrow\infty$ in \eqref{contr}. Therefore, for $\varepsilon<\varepsilon_0$, we obtain the desired upper bound,
$$T\leq C_4\varepsilon^{-\frac{2p(p-1)}{\gamma(p,n+\mu)}}$$
with
$$C_4:=\left(\frac{e^{(S_p(\infty)+\alpha\log2)+1}}{C_2}\right)^{2(p-1)/\gamma(p,n+\mu)}.$$
This completes our proof of main theorem.
\end{proof}

{\bf Acknowledgment}:
The authors would like to thank Professor Takamura for his helpful comments and suggestions which lead to remove the restriction of $p>p_F(n+\frac\mu2)$ in the previous manuscript of arXiv:1709.00866. This work was done when the the first author visited Mathematical Department, Hokkaido University. He would like thank Professor Jimbo's kindly invitation and great support. The second author is partially supported by NSFC No. 11501511 and Zhejiang Provincial Nature Science Foundation of China under Grant No. LQ15A010012.

\newpage
\bibliographystyle{elsarticle-num}
\bibliography{References}

\begin{thebibliography}{99}

\bibitem{Dabbi15JDE} D'Abbicco M, Lucente S, Reissig M. A shift in the Strauss exponent for semilinear wave equations with a not effective damping[J]. Journal of Differential Equations, 2015, 259(10): 5040-5073.
\bibitem{Erdelyi} Erdelyi, A., Magnus, W., Oberhettinger, F., Tricomi, F.G.: Higher Transcendental Functions, vol. 2.
McGraw-Hill, New York (1953)
\bibitem{Tri1} He D, Witt I, Yin H. On the global solution problem for semilinear generalized Tricomi equations, I[J]. Calculus of Variations and Partial Differential Equations, 2017, 56(2): 21.
\bibitem{Ikedapre} Ikeda M, Sobajima M. Life-span of solutions to semilinear wave equation with time-dependent critical damping for specially localized initial data[J]. arXiv preprint arXiv:1709.04406, 2017.
\bibitem{Lai17} Lai N A, Takamura H, Wakasa K. Blow-up for semilinear wave equations with the scale invariant damping and super Fujita exponent[J]. arXiv preprint arXiv:1701.03232, 2017.
\bibitem{Lai-scattering} Lai N A, Takamura H. Blow-up for semilinear damped wave equations with sub-Strauss exponent in the scattering case[J]. arXiv preprint arXiv:1707.09583, 2017.
\bibitem{Takam} Takamura H. Improved Kato¡¯s lemma on ordinary differential inequality and its application to semilinear wave equations[J]. Nonlinear Analysis: Theory, Methods \& Applications, 2015, 125: 227-240.
\bibitem{Wakasugi14}Wakasugi Y. Critical exponent for the semilinear wave equation with scale invariant damping[M]//Fourier Analysis. Birkh\"{a}user, Cham, 2014: 375-390.
\bibitem{WakasugiT} Wakasugi Y. On the diffusive structure for the damped wave equation with variable coefficients[D]. PhD Thesis, Osaka University, 2014.
\bibitem{Wakasa16} Wakasa K. The lifespan of solution to the semilinear damped wave equations in one space dimemsion[J]. Communications on Pure \& Applied Analysis, 2016, 15(4).
\bibitem{Wirt04} Wirth J. Solution representations for a wave equation with weak dissipation[J]. Mathematical methods in the applied sciences, 2004, 27(1): 101-124.
\bibitem{Zhang} Yordanov B T, Zhang Q S. Finite time blow up for critical wave equations in high dimensions[J]. Journal of Functional Analysis, 2006, 231(2): 361-374.
\bibitem{Zhou14} Zhou Y, Han W. Life-span of solutions to critical semilinear wave equations[J]. Communications in Partial Differential Equations, 2014, 39(3): 439-451.
\end{thebibliography}

\end{document}